# Quantifying extinction probabilities of endangered species for phylogenetic conservation prioritization may not be as sensitive as might be feared

Alain Billionnet

*Laboratoire CEDRIC, École Nationale Supérieure d'Informatique pour l'Industrie et l'Entreprise,*
*1, square de la Résistance, 91025 Évry cedex, France; E-mail: Alain.Billionnet@ensiie.fr.*

*Abstract*. In this study we are concerned with the general problem of choosing from a set of endangered species $T$ a subset $S$ of $k$ species to protect as a priority. Here, the interest to protect the species of $S$ is assessed by the resulting expected phylogenetic diversity (ePD) of the set $T$, a widely used criterion for measuring the expected amount of evolutionary history associated with $T$. We consider that the survival of the protected species is assured and, on the contrary, that there is a risk of extinction for the unprotected species. The problem is easy to solve by a greedy type method if the extinction probabilities of the unprotected species are known but these probabilities are generally not easy to quantify. We show in this note that the choice of the precise values attributed to the extinction probabilities – provided it respects the rank of imperilment of each species – is not as decisive as might be feared for the considered problem. The values of these probabilities have a clear impact on the selection of the species to be protected but a little impact on the resulting ePD. More precisely, if $T_1$ and $T_2$ are the two optimal subsets of species corresponding to two scenarios (two different sets of probabilities) the ePDs of $T_1$ and $T_2$, calculated with the probabilities of the first scenario – or with the probabilities of the second scenario – are not very different.

1.  **Introduction**

An increasing number of species are threatened with extinction (Stork 2010). The causes for this growing threat include habitat loss, degradation, and fragmentation, agriculture, human over-population, deforestation, poaching and hunting, invasive species, and climate change (Collen et al. 2013). The loss of biodiversity that would result from a massive extinction would have disastrous consequences (Cardinale et al. 2012). The international community agrees that the protection of threatened and endangered species is a priority issue but the resources available to implement these protections are limited (Convention on Biological Diversity 2011). It is therefore crucial to carefully select the species on which the protection efforts will be concentrated in order to preserve biodiversity as much as possible. There are many ways to assess biodiversity (see, for example, (Moreno et al. 2017) and also the very extensive book by MacLaurin and Sterelny (2008) on this subject). In this paper, the measurement of biodiversity is based on the notion of phylogenetic diversity. We consider that the survival of a protected species is assured - which is not the case for an unprotected species - and we consider the classic conservation problem which consists in determining, among a set of species $T$, the "best" subset of species to protect within a limited budget. More precisely, we consider the problem of determining, among a set of species $T$, $k$ species to be protected so as to maximize the resulting expected phylogenetic diversity of $T$, $k$ being a fixed value. The phylogenetic diversity is a widely used criterion for expressing the amount of biodiversity (evolutionary history) of a set of species captured by a phylogeny. The reader can refer to (Maclaurin & Sterelny 2008, Chap.7 and Pellens & Grandcolas 2016) for a very comprehensive presentation of the role of phylogenetic diversity, in large sense, in the field of biological conservation. The expected phylogenetic diversity requires the precise knowledge of the extinction probabilities of the species concerned and these probabilities are generally not easy to quantify (Mace et al. 2008 ; Mooers et al. 2008). In this note, we use the measure of phylogenetic diversity created by Faith (Faith 1992a, 1992b) and we study the influence, on the solutions of the considered conservation problem, of the values attributed to the extinction probabilities.

2. **The expected phylogenetic diversity**

A rooted phylogenetic tree can be considered as a directed tree, $H = (N, A, T, \lambda)$, where $N$ is the set of nodes, $A$, the set of arcs – or directed edges or branches –, $T$, the set of leaves representing the species, and $\lambda$, the vector of the length of each branch. The branches of $A$ are directed from the root of $H$ to its leaves. The nodes of the phylogenetic tree that are neither root nor leaves are called internal nodes. The internal nodes have one predecessor and at least two successors. For all $a \in A$, we denote by $L_a$ the set of species for which there is a path from the terminal node of $a$ to these species. As mentioned in Section 1, we use the measure of phylogenetic diversity created by Faith: The phylogenetic diversity of the set $T$ of species associated with the tree $H$ is equal, by definition, to $\sum_{a \in A} \lambda_a$. Intuitively, it represents the total amount of evolutionary history embodied in a set of species, since the time of the most recent common ancestor of the set (Chao et al. 2010). In this paper, we refer to PD to indicate this measure. The reader is referred to (Faith 2013, 2016) and also (Collen 2015; Isaac et al. 2007) for a comprehensive discussion about the use of PD in measuring biodiversity. Mimouni et al. (2016) point out that PD is an established index with a large enough body of literature regarding its computational and combinatorial aspects. For example, Hartmann & Steel (2007) review algorithmic, mathematical, and stochastic results concerning PD. As mentioned in the introduction, the book of Maclaurin and Sterelny (2008) and the book edited by Pellens and Grandcolas (2016) deal in detail with the importance of phylogenetic diversity, and especially Faith's phylogenetic diversity, for biodiversity conservation.

Suppose that an extinction probability, $p_i$, is associated with each species, $i$, of $T$. The probability that the information associated with an arc $a$ of the tree is retained is equal to the probability that at least one of the species "reachable" from $a$ – i.e. belonging to $L_a$ – survives. This way, we can define (from the Faith's PD) the expected phylogenetic diversity (ePD) of $T$. It is equal to the sum on all arcs of $T$ of the probability that the information associated with this arc is conserved, multiplied by the length of this arc (Witting & Loeschcke 1995). More precisely, the ePD of $T$ is equal to $\sum_{a \in A} \lambda_a (1 - P_a)$ where $P_a = \prod_{i \in L_a} p_i$. We suppose here that the extinction probabilities of the species are independent. The ePD is recognized as a relevant criterion for evaluating the biodiversity associated with a set of species. It accounts for complementarity among species (Faith 2015).

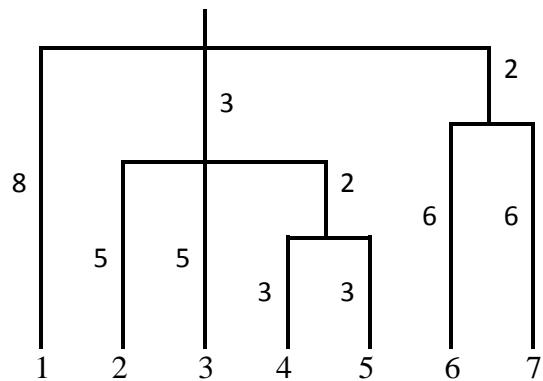

**Figure 1.** A hypothetical ultrametric phylogenetic tree associated with 7 species; it has 10 branches. The length of each branch is indicated next to the branch.

As an example, consider the hypothetical ultrametric phylogenetic tree of Figure 1 and assume that the 7 species considered all have the same extinction probability. Let us examine the impact of the following 2 events (of same probability): the extinction of the 3 species 1, 6 and 7, the 4 other species remaining alive, and the extinction of the 3 species 2, 3 and 4, the 4 other species remaining alive. The first event would cause the loss of 22 units while the second would cause the loss of only 13 units.

The use of ePD requires the knowledge of the extinction probabilities, $p_i$, for each species *i*, and these probabilities are generally not easy to quantify (Mooers et al. 2008). We consider that the studied species are divided into several groups more or less threatened and that an extinction probability is associated with each group. It is possible, for example, to consider the 5 IUCN Red List categories: Least Concern, Near Threatened, Vulnerable, Endangered and Critically Endangered. Note that ranking species in these categories is a difficult task (Mace et al. 2008). The extinction probability of the species *i*, $p_i$, is then equal to the extinction probability assigned to the group of threat to which it belongs. Redding & Mooers (2006) and Isaac et al. (2007) have proposed ways to transform the Red List categories to extinction probabilities. Mooers et al. (2008) have shown that the choice of transformation of Red List categories to different extinction probabilities can affect the resulting species ranking for the metrics EDGE (Isaac et al. 2007) and HEDGE (Steel et al. 2007; Mooers et al. 2008; Forest et al. 2015; Nunes et al. 2015). In this study, we examine the impact of the extinction probability value attributed to each group (and thus to each species) on the optimal protection strategy associated with the classical conservation problem outlined above and precisely defined below.

## 3. Determining the best set of species to protect

The problem that we consider here can be stated as follows: given a phylogenetic tree $H = (N, A, T, \lambda)$ and extinction probabilities, $p_i$, for each species, *i*, of *T*, determine a set of *k* species of *T* whose protection maximizes the ePD of *T*. Remember that we assume that the protection of a species ensures its survival. In other words, its probability of extinction decreases from $p_i$ to 0. This basic problem is a particular case of the Noah's Ark problem which was first defined by Weitzman (1998). It is shown in (Hartman & Steel 2006) and (Pardi 2009) that the considered conservation problem can be solved easily by a greedy algorithm. The algorithm begins with an empty set, *S*, and sequentially adds to *S* the species, *i*, which maximizes the ePD of the set $S \cup \{i\}$ until the cardinal of *S* is *k*.

*Link with the HEDGE metric*: The well-known HEDGE metric defines species priority ranks for a fixed tree and fixed extinction probabilities. The HEDGE score of the species $i$ is equal to the increase of the ePD of $T$ obtained by protecting this species with respect to the ePD of $T$ obtained in protecting no species (Figure 2). One way of approaching the considered problem - determine the species of $T$ to be protected so as to maximize the resulting ePD of $T$- consists in calculating the HEDGE score of each species and then retaining the $k$ species with the highest scores. It is to be noted, however, that the HEDGE score does not take into account the protection actions actually taken. Indeed, if the species that has the highest HEDGE score is effectively protected (its extinction probability becomes equal to 0) the HEDGE scores of some other species can decrease significantly. To overcome this drawback, a natural extension of the HEDGE method, called I-HEDGE, has been proposed by Jensen et al. (2016). This method consists in gradually calculating the HEDGE score of each species taking into account at each stage the fact that a new species becomes protected. The new species that becomes protected is the one with the highest score, this one being calculated taking into account the species already chosen to be protected. This method, which therefore integrates the notion of phylogenetic complementarity, corresponds exactly to the greedy algorithm proposed by Pardi (2009) and Hartman & Steel (2006), and mentioned above. Figure 2 illustrates the computation of the EDGE and I-HEDGE scores.

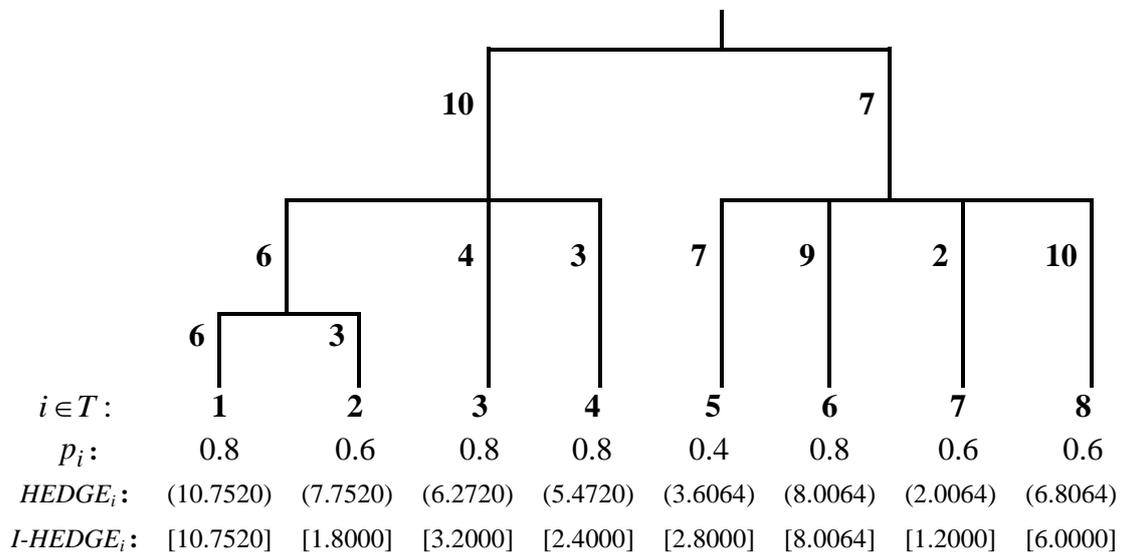

**Figure 2.** A hypothetical phylogenetic tree with 8 species. The lengths of the branches are indicated next to each branch. The extinction probabilities of each unprotected species ($p_i$) are shown under each species. The HEDGE score for each species (*HEDGE$_i$*) is shown in brackets under the corresponding species. The I-HEDGE score for each species (*I-HEDGE$_i$*) is shown in square brackets under the corresponding species. According to these last values, the 3 species whose protection maximizes the ePD of $T$ are 1, 6 and 8, and the corresponding ePD is equal to 55.6 (30.8416+10.7520+8.0064+6.0000).

Consider the detail of the calculation of the I-HEDGE scores on the tree in Figure 2. If no species is protected, the ePD of *T* is equal to 30.8416. The species whose protection increases the most the ePD of *T* is species 1. The protection of this species increases the ePD of *T* from 30.8416 to 30.8416+10.7520=41.5936. Let us definitely set the value of the survival probability of species 1 to 1. The species (different from 1) whose protection increases the most the ePD of *T* is species 6. The protection of this new species increases the ePD of *T* from 41.5936 to 41.5936+8.0064=49.6000. Let us definitively set the survival probability of species 6 to 1. The species (different from 1 and 6), the protection of which increases the most the ePD of *T* is species 8. The protection of this new species increases the ePD of *T* from 49.6000 to 49.6000+6.0000=55.6000.

## 4.    Main result

We present in this section experimental results which show that the choice of the values attributed to the extinction probabilities associated with the 5 groups of species considered does not have an immense influence on the solution of the problem under consideration although this choice may significantly affect the set of species to be protected as a priority. Consider the two scenarios defined by the probabilities $\pi_g^\omega$ ($g=1,...,5; \omega=1,2$) which designate the extinction probability of the species of the group *g* in the scenario $\omega$. We denote by $S^{*\omega}$ ($\omega=1,2$) the optimal set of species to be protected in scenario $\omega$. As we have seen above, these sets can be determined by the greedy algorithm. For each subset *S* of *T* and for each scenario $\omega$, we denote by $ePD^\omega(T,S)$ the ePD of *T* generated by the protection of the species of *S* and calculated with the probabilities of the scenario $\omega$.

We are interested first of all in the case where one decides to protect the species of $S^{*1}$. We then compute the resulting ePD of *T* with the two sets of extinction probabilities, that is to say $ePD^1(T,S^{*1})$ and $ePD^2(T,S^{*1})$. We then consider the case where one decides to protect the species of $S^{*2}$ and one calculates the resulting ePD of *T* with the two sets of extinction probabilities, that is to say $ePD^1(T,S^{*2})$ and $ePD^2(T,S^{*2})$. We then calculate the two relative gaps denoted by $gap^1$ and $gap^2$:

$$gap^1 = [ePD^1(T,S^{*1}) - ePD^1(T,S^{*2})]/ePD^1(T,S^{*1})$$

$$gap^2 = [ePD^2(T,S^{*2}) - ePD^2(T,S^{*1})]/ePD^2(T,S^{*2})$$

The first of these gaps provides the relative error made (on the maximum ePD) in considering that the real scenario is the scenario 2 when in reality it is the scenario 1. The second gap provides the relative error made in considering that the real scenario is the scenario 1 when in reality it is the scenario 2.

We first calculated $gap^1$ and $gap^2$ on 10,000 randomly generated instances corresponding to non ultrametric trees and also on 10,000 randomly generated instances corresponding to ultrametric trees. The description of all these instances and the corresponding detailed results are presented in Section 5.

For the non ultrametric trees, the largest of the gaps obtained (out of 20,000) is equal to 2.63% (Table 5). Therefore, we can say that, at least under these experimental conditions, the possible error committed – on the maximum ePD values – by the choice of transformation of categories to different extinction probabilities is relatively limited. However, it should be noted that the set of *k* species, optimal in scenario 1, may be quite different from the set of *k* species optimal in scenario 2. We measured their dissimilarity by the cardinal of their symmetric difference divided by the cardinal of their union. This criterion therefore varies from 0 to 1. For example, the dissimilarity between the two optimal sets corresponding to the largest gap is equal to 0.70 which corresponds to two very different sets (Table 5). Note that Mooers et al. (2008) had already observed that the choice of transformation of Red List categories to different extinction probabilities can affect significantly the resulting species ranking for the two well-known metrics EDGE and HEDGE.

For the ultrametric trees, the largest of the gaps obtained (out of 20,000) is slightly higher than the one obtained for non-ultrametric trees since it is equal to 4.13% (Table 5). Here too, the dissimilarity between the two optimal sets corresponding to this gap is important. It is equal to 0.59 (Table 5).

We then calculated $gap^1$ and $gap^2$ on 10,000 instances generated from the phylogenetic tree associated with the 52 Madagascar lemurs (species and subspecies) and studied by Pardi (2009). In this case, the phylogenetic tree considered (structure, length of the branches and category to which each of the 52 lemurs belongs) is definitively fixed. It includes 97 nodes in total. Each of the 10,000 instances is generated by associating to the 52 lemurs 2 different scenarios of extinction probabilities. The description of these instances and the corresponding detailed results are presented in Section 5. As in experiments on randomly generated phylogenetic trees we found that 1) the possible error committed – on the maximum ePD values – by the choice of transformation of categories to different extinction probabilities is relatively limited (0.89%) and 2) the two optimal sets corresponding to the largest error are very different (dissimilarity equal to 0.42).

**5. Description of the instances and detailed results**

*Randomly generated non ultrametric phylogenetic trees*

The 10,000 instances considered were generated in the following way:

*Construction of the hypothetical phylogenetic trees*: First, the number of non-leaf nodes is uniformly and randomly chosen between 50 and 1,000. Then, for each of these nodes, the number of children is uniformly and randomly generated between 2 and $d_{max}$, and $d_{max}$ is itself drawn at random uniformly between 2 and 4.

*Lengths of the branches*: The length of each branch is uniformly and randomly generated between 1 and $\lambda_{max}$, and $\lambda_{max}$ is itself drawn at random uniformly between 5 and 20.

*Assignment of extinction probabilities*: The set of species is divided into five categories and the category to which each species belongs is drawn at random in two different ways. In the first case, this group is drawn at random uniformly in the set {1,2,3,4,5}. The numbers of species in each group are then comparable. In the second case, the group is drawn at random to approximately obtain the following distribution of the species in the five groups, as in (Mooers et al. 2008): 2%, 4%, 9%, 9%, 76%. Considering the case 1 or the case 2 is chosen at random in an equiprobable manner. To define the extinction probabilities of an unprotected species we were inspired by the description of the 5 IUCN Red List categories (IUCN 2012): For each tree considered, two sets of extinction probabilities (corresponding to two hypothetical scenarios) are drawn at random uniformly in the intervals given in Table 1.

**Table 1.** Randomly generated trees: possible values of the extinction probability of each category.

| Category | Interval in which the extinction probability is chosen at random |
|---|---|
| 1 | [0.50,1.00] |
| 2 | [0.20,0.50] |
| 3 | [0.10,0.20] |
| 4 | [0.05,0.10] |
| 5 | [0.00,0.05] |

*Maximum number of species that can be protected*: Finally, the maximum number of species that can be protected is equal to the greatest integral number less than or equal to $\rho$ multiplied by the number of species, and $\rho$ is drawn at random uniformly in the set {0.1, 0.2, 0.3, 0.4, 0.5}.

We also conducted other experiments slightly different from the previous ones by modifying the way extinction probabilities are attributed. In these experiments, the extinction probabilities of the species may vary within a category. More precisely, the extinction probability of each species is drawn at random uniformly in the range corresponding to its category. For example, if species 27 belongs to category 2, its extinction probability is drawn at random uniformly in the interval [0.20, 0.50]. We again generated 10,000 instances and we obtained results comparable to those we had previously obtained: The mean of the 20,000 relative gaps is equal to 0.19%, the corresponding standard deviation is equal to 0.13%, and the largest relative gap is equal to 1.21%. The instance corresponding to this largest gap concerns 67 species divided into five categories respectively formed of 16, 12, 10, 15 and 14 species. In this instance, the maximum number of species that can be protected is equal to 6 and the dissimilarity of the two optimal sets obtained is equal to 0.50.

*Randomly generated ultrametric phylogenetic trees*

The 10,000 instances considered were generated in exactly the same way as for non-ultrametric trees, and then the lengths of the branches linked to a species in each resulting tree are modified to obtain an ultrametric tree. Let *L* be the length of the longest path connecting the root to a species and let $l_i$ be the length of the path connecting the root to the species *i*. For each species *i*, the length of the arc connected to this species is modified by adding to it the value *L-$l_i$*.

*Phylogenetic tree associated with Madagascar lemurs (Pardi 2009)*

The 10,000 instances considered were generated in the following way:
*Assignment of extinction probabilities*: The set of species is divided into five categories. The category to which each species belongs is given by Pardi (2009). For each instance considered, two sets of extinction probabilities (corresponding to two hypothetical scenarios) are drawn at random uniformly in the intervals given in Table 2. These intervals represent a variation around the values of the extinction probabilities proposed by Pardi (2009).

**Table 2.** Madagascar lemurs: possible values of the extinction probability of each category.

| Category | Interval in which the extinction probability is chosen at random |
|---|---|
| 1 | [0.90,1.00] |
| 2 | [0.65,0.85] |
| 3 | [0.40,0.60] |
| 4 | [0.15,0.35] |
| 5 | [0.00,0.10] |

*Maximum number of species that can be protected*: As in experiments on randomly generated phylogenetic trees, the maximum number of species that can be protected is equal to the greatest integral number less than or equal to $\rho$ multiplied by the number of species (52), and $\rho$ is drawn at random uniformly in the set $\{0.1, 0.2, 0.3, 0.4, 0.5\}$.

For the three types of phylogenetic trees considered (non-ultrametric, ultrametric, Madagascar lemur), statistics on the 10,000 generated instances, characteristics of the instance corresponding to the largest gap and solutions of the conservation problem considered for this instance and for the two associated scenarios are presented in Tables 3-5.

With regard to the phylogenetic tree associated with Madagascar lemurs we also measured the impact that possible errors on the lengths of the branches could have on the solution of the considered conservation problem: determine $k$ species to protect in order to maximize the resulting ePD. For this we have definitively assigned to each species the survival probability given by Pardi (2009) and we have randomly generated two sets of values for the lengths of the branches, thus obtaining two scenarios. The value of each branch is uniformly and randomly generated in the interval $[0.75\lambda_a, 1.25\lambda_a]$ where $\lambda_a$ is the value of the branch given by Pardi. We therefore consider that a significant error ($\pm 25\%$) can be made on the lengths of the branches. We have thus built 10,000 different instances. For each instance, we have determined the 2 optimal sets $S^{*1}$ and $S^{*2}$ corresponding to the two scenarios and we have calculated the two relative gaps, $gap^1$ and $gap^2$, defined in Section 4. Here, $ePD^\omega(T,S)$ is the ePD of $T$ generated by the protection of the species of $S$ and calculated with the branch lengths of the scenario $\omega$ ($\omega$=1 or 2). The maximum of these 2 gaps over the 10,000 instances generated, that is to say the maximum error that can be made by choosing the wrong scenario, is equal to 1.68%. This error is larger than that obtained by considering uncertainties on the survival probabilities. Moreover the two optimal sets of species corresponding to the largest gap are very different since their dissimilarity is equal to 0.64. If we admit that an error of $\pm 50\%$ can affect the length of the branches, the maximum relative gap goes to 5.28% and the dissimilarity to 0.57.

**Table 3.** Statistics on the 10,000 randomly generated instances.

|  | Non-ultrametric trees | Ultrametric trees | Madagascar lemurs tree |
|---|---|---|---|
| Maximum number of species in the 10,000 generated trees | 2058 | 2058 | - |
| Minimum number of species in the 10,000 generated trees | 74 | 74 | - |
| Mean of the 20,000 relative gaps (%) | 0.08 | 0.07 | 0.08 |
| Standard deviation of the 20,000 relative gaps (%) | 0.14 | 0.20 | 0.12 |
| Maximum of the 20,000 relative gaps (%) | 2.63 | 4.13 | 0.89 |
| Dissimilarity maximum of the 10,000 pairs of optimal sets | 0.70 | 0.68 | 0.42 |

**Table 4.** Characteristics of the instances corresponding to the largest relative gap.

|  | Non-ultrametric trees | Ultrametric trees | Madagascar lemurs tree |
|---|---|---|---|
| Number of non-leaf nodes | 727 | 147 | - |
| Number of species | 1420 | 290 | - |
| $d_{max}$ | 4 | 4 | - |
| $\lambda_{max}$ | 7 | 9 | - |
| Length of the longest path from the root to a species | 72 | - | - |
| Length of the shortest path from the root to a species | 16 | - | - |
| Number of species in each group | 297, 293, 267, 274, 289 | 55, 64, 50, 58, 63 | - |
| Scenario 1: Extinction probabilities by category: | 0.86, 0.24, 0.14, 0.07, 0.01 | 0.51, 0.47, 0.12, 0.06, 0.01 | 0.99, 0.76, 0.48, 0.32, 0.05 |
| Scenario 2: Extinction probabilities by category: | 0.50, 0.46, 0.13, 0.07, 0.03 | 0.99, 0.22, 0.18, 0.07, 0.03 | 0.95, 0.67, 0.59, 0.16, 0.01 |
| $\rho$ | 0.1 (142 species). | 0.2 (58 species) | 0.3 (15 species). |
| Scenario 1: ePD obtained by protecting no species | 11003.39 | 8262.40 | 94.07 |
| Scenario 2: ePD obtained by protecting no species | 11237.61 | 7748.15 | 98.70 |
| ePD obtained by protecting all the species | 13485.00 | 10520.00 | 143.18 |

**Table 5.** Results associated with the solution of the problem
in the two scenarios for the instances described in Table 4.

|  | Non-ultrametric trees | Ultrametric trees | Madagascar lemurs tree |
|---|---|---|---|
| Scenario 1: Maximum ePD | 12176.37 | 9407.09 | 131.21 |
| Scenario 1: Distribution by category of the species to be protected | category 1: 142 | category 1: 31<br>category 2: 27 | category 1: 5<br>category 2: 6<br>category 3: 2<br>category 4: 2 |
| ePD obtained in scenario 1 by protecting the species of the set optimal in scenario 2: | 11856.61 | 9265.58 | 130.87 |
| $gap^1$ (%) | 2.63 | 1.50 | 0.26 |
| Scenario 2: Maximum ePD | 12001.33 | 9605.16 | 132.97 |
| Scenario 2: Distribution by category of the species to be protected | category 1: 66<br>category 2: 76 | category 1: 55<br>category 2: 3 | category 1: 4<br>category 2: 6<br>category 3: 5 |
| ePD obtained in scenario 2 by protecting the species of the set optimal in scenario 1 | 11914.76 | 9208.83 | 131.79 |
| $gap^2$ (%) | 0.72 | 4.13 | 0.89 |
| Dissimilarity of the two optimal sets | 0.70 | 0.59 | 0.42 |

## 6. Conclusions and perspectives

To determine which of a given set of endangered species should be protected as priorities, the criterion of the expected phylogenetic diversity (ePD) is recognized as a relevant criterion. However, its use requires explicit knowledge of the extinction probabilities of the species considered when they are protected and also when they are not protected. We assume here that the extinction probability of the protected species is zero. We were interested in the influence of the conversion of the extinction risk of a species to an extinction probability of this species. The extinction risk of a species is given by the fact that it belongs to a certain category (such as those of the IUCN Red List). We showed experimentally, on randomly generated ultrametric and non-ultrametric phylogenetic trees as well as on a real tree, that the choice of the extinction probabilities values is not as important as might be feared. More precisely, if we consider 2 scenarios, corresponding to two different conversions, the ePD generated by the protection of a set of $k$ species, optimal in the scenario 1 and calculated with the probabilities of this scenario is not very distant from the ePD generated by the protection of a set of $k$ species, optimal in the scenario 2 and calculated with the probabilities of the scenario 1. In other words, the error committed by retaining the scenario 2 when the real scenario would be the scenario 1 is not very important, at least in our experimental conditions. The terms "scenario 1" and "scenario 2" can be exchanged in the preceding two sentences. It should be noted, however, that the set of *k* species, optimal in the scenario 1, may be quite different from the set of *k* species optimal in the scenario 2.

We considered in this study that the survival probabilities of the protected species were equal to 1. It would be interesting to study other values for these probabilities and to measure the impact of the values retained (if the species is protected and if it is not protected) on the optimal ePD. It would also be interesting to study the impact of the uncertainties that inevitably exist on the structure of the phylogenetic trees considered and the values attributed to their branches. In section 5, we have given some preliminary results that show that the consequences of uncertainties on branch lengths may not be negligible.

**References**

Cardinale B.J., Duffy J.E., Gonzalez A., Hooper D.U., Perrings C., Venail P., Narwani A., Mace G.M., Tilman D., Wardle D.A., Kinzig A.P., Daily G.C., Loreau M., Grace J.B., Larigauderie A., Srivastava D.S., Naeem S. 2012. Biodiversity loss and its impact on humanity. Nature, 486, 59-67.

Chao A., Chiu C.-H., Jost L. 2010. Phylogenetic diversity measures based on Hill numbers. Phil. Trans. R. Soc. B 365, 3599–3609.

Collen B. 2015. Conservation prioritization in the context of uncertainty. Animal Conservation, 18, 315–317.

Collen B., Pettorelli N., Baillie J.E.M., Durant S.M. (Eds) 2013. Biodiversity monitoring and conservation: Bridging the gap between global commitment and local Action. Wiley-Blackwell.

Convention on Biological Diversity (2011). Strategic plan for biodiversity 2011–2020, including Aichi biodiversity targets. Available at: http://www.cbd.int/sp/.

Faith D.P. 1992a. Conservation evaluation and phylogenetic diversity. Biol. Cons., 61, 1-10.

Faith D.P. 1992b. Systematics and conservation: on predicting the feature diversity of subsets of taxa. Cladistics, 8, 361–373.

Faith D.P. 2013. Biodiversity and evolutionary history: useful extensions of the PD phylogenetic diversity assessment framework. Annals of the New-York Academy of Sciences, 1289, 69–89.

Faith D.P. 2015 Phylogenetic diversity, functional trait diversity and extinction: avoiding tipping points and worst-case losses. Phil. Trans. R. Soc. B370: 20140011.

Faith DP. 2016. The PD phylogenetic diversity framework: linking evolutionary history to feature diversity for biodiversity conservation. In: Pellens R., Grandcolas P. (eds), Biodiversity Conservation and Phylogenetic Systematics, Topics in Biodiversity and Conservation 14, Springer.

Forest F., Crandall K.A., Chase M.W., Faith D.P. 2015. Phylogeny, extinction and conservation: embracing uncertainties in a time of urgency. Phil. Trans. R. Soc. B 370:20140002.

Hartmann K., Steel M. 2006. Maximizing phylogenetic diversity in biodiversity conservation: greedy solutions to the Noah's Ark Problem. Syst. Biol., 55, 644–651.

Hartmann K., Steel M. 2007. Phylogenetic diversity: From combinatorics to ecology. In: Gascuel O. & Steel M. (Eds), Reconstructing evolution: new mathematical and computational advances. Oxford University Press, Oxford.

Isaac N.J.B., Turvey S.T., Collen B., Waterman C., Baillie J.E.M. 2007. Mammals on the EDGE : Conservation priorities based on threat and phylogeny. PloS ONE 2(3): e296.

IUCN. (2012). IUCN Red List Categories and Criteria: Version 3.1. Second edition. Gland, Switzerland and Cambridge, UK: IUCN. iv + 32pp.

Jensen E.L., Mooers A.O., Caccone A., Russello M.A. 2016. I-HEDGE: determining the optimum complementary sets of taxa for conservation using evolutionary isolation. PeerJ, 4: e2350.

Mace G.M., Collar N.J., Gaston K.J., Hilton-Taylor C., Akçakaya H.R., Leader-Williams N., Milner-Gulland E.J., Stuart S.N. 2008. Quantification of extinction risk: IUCN's system for classifying threatened species. Conservation Biology, 22, 1424-1442.

MacLaurin J., Sterelny K. (2008). What Is Biodiversity? The University of Chicago Press, Chicago, 217 p.

Mimouni E., Beisner B.E, Pinel-Alloul B.. 2016. Phylogenetic diversity and its conservation in the presence of phylogenetic uncertainty: a case study of cladoceran communities in urban waterbodies. Biodiversity and Conservation, 25, 2113-2136.

Mooers A.O., Faith D. P., Maddison W.P. 2008. Converting endangered species categories to probabilities of extinction for phylogenetic conservation prioritization. PLoS One 3:e3700.

Moreno C.E. et al. 2017. Measuring biodiversity in the Anthropocene: a simple guide to helpful methods. Biodiversity and Conservation, to appear.

Nunes L.A., Turvey S.T., Rosindell J. 2015. The price of conserving avian phylogenetic diversity: a global prioritization approach. Phil. Trans. R. Soc. B 370.

Pardi F. 2009. Algorithms on phylogenetic trees. PhD thesis, University of Cambridge.

Pellens R., Grandcolas P. (Eds). 2016. Biodiversity conservation and phylogenetic systematics: preserving our evolutionary heritage in an extinction crisis. Springer International Publishing, Springer Open, 390 p.

Redding R.W., Mooers A.O. 2006. Incorporating phylogenetic measures into conservation prioritization. Cons. Biol., 20, 1670–1678.

Steel M., Mimoto A., Mooers A.O. 2007. Hedging our bets: The expected contribution of species to future phylogenetic diversity. Evol. Bioinform., 3, 237–244.

Stork N.E. 2010. Re-assessing current extinction rates. Biodivers. Conserv., 19, 357–371.

Weitzman M.L. 1998. The Noah's ark problem. Econometrica 66,1279–1298.

Witting L., Loeschcke V. 1995. The optimization of biodiversity conservation. Biol. Conserv., 71, 205–207.